\documentclass[12pt]{article}
\newtheorem{theor}{Theorem}
\newtheorem{prop}{Proposition}

\begin{document}

\title{Non-rigidity degrees of root lattices and their duals}

\author{
Michel Deza\footnote{Michel.Deza@ens.fr}
\\
Ecole Normale Sup\'erieure, Paris, and ISM, Tokyo
\and
Viacheslav Grishukhin\footnote{grishuhn@cemi.rssi.ru} \\
CEMI, Russian Academy of Sciences, Moscow}

\date{}

\maketitle

\begin{abstract}
Non-rigidity degree of a lattice $L$, nrd$L$, is dimension of 
the L-type domain to which $L$ belongs. We complete here the table 
of nrd's of all root lattices and their duals; namely, the hardest
remaining case of $D_n^*$, and the case of $E_7^*$ are decided. 

We describe explicitly the $L$-type domain ${\cal D}(D_n^*)$, $n \ge 4$. 
For $n$ odd, it is a non-simplicial polyhedral open cone of dimension 
$n$. For $n$ even, it is one-dimensional, i.e. for even $n$, $D_n^*$ 
is an edge form.
\end{abstract}

{\em Mathematics Subject Classification}. Primary 11H06, 11H55;
Secondary 52B22, 05B45.

{\em Key words}. Root lattices, Delauney
polytopes, rank.

\section{Introduction}
Voronoi \cite{Vo} defined the partition of the cone ${\cal P}_n$
of positive semidefinite $n$-ary quadratic forms into $L$-type
domains, that we call here {\em $L$-domains}. Forms of the
same $L$-domain correspond to lattices that determine affinely equivalent
Voronoi partitions of ${\bf R}^n$, i.e. partition into Voronoi polytopes.
The partition of ${\bf R}^n$, which is (both, combinatorially and affinely)
dual to a Voronoi partition, is called Delaunay partition and consists of
Delaunay polytopes.  In other words,
two lattices have the same $L$-type if and only if the face posets
of their Voronoi polytopes are isomorphic. (The {\em face poset} of
a polytope $P$ is the set of all faces of $P$ of all dimensions
ordered by inclusion.)

If the $L$-type of a lattice changes, then either some Delaunay
polytopes are glued into a new Delaunay polytope, or a Delaunay
polytope is partitioned into several new Delaunay polytopes.
Recall that the center of a Delaunay polytope is a vertex of a
Voronoi polytope. Hence if the $L$-type of a lattice changes, then
for each Voronoi polytope either some vertices are glued into one
vertex, or a vertex splits into several new vertices.

Voronoi proved that each $L$-domain is an open polyhedral cone of
the dimension $k$, $1 \le k \le N$, where $N=\frac{n(n+1)}{2}$ is
the dimension of ${\cal P}_n$, i.e. dimension of the space of
coefficients of $f \in {\cal P}_n$. An $L$-type having an
$N$-dimensional domain, is called {\em generic}.
Otherwise, $L$-type is called {\em special}.

In \cite{BG}, a notion of {\em non-rigidity degree} of a form $f$
and the corresponding lattice $L(f)$ is introduced. It is denoted by
nrd$f$ and is equal to dimension of the $L$-domain containing
$f$. It is shown in \cite{BG} that nrd$f$ is equal to corank of a
system of equalities connecting the norms of minimal vectors of cosets
$2L$ in $L$. In fact, nrd is the number of degrees of freedom, such 
that the Delaunay partition (more precisely, its star) has, when one 
deforms it affinely, so that the result remains an Delaunay partition. 
Clearly, the maximal nrd is $n+1 \choose 2$ and it is realized only by 
simplicial Delaunay partitions.

So, $1 \le {\rm nrd}f \le N$, and nrd$f=N$ if $f$ belongs to a generic
$L$-domain. A form $f$ and the corresponding lattice $L(f)$ are called 
{\em rigid} if nrd$f=1$. This name was used because any, distinct from 
a scaling, affine transformation of a rigid lattice changes its $L$-type. 
(Sometimes, a rigid form is called an {\em edge form}, since it lies on 
an extreme ray of the closure of an $L$-domain.) Clearly, any 
1-dimensional lattice is trivially rigid. In \cite{BG}, were given 7 
examples of rigid lattices of dimension 5 and shown that $D_4$ is unique 
such a lattice of dimension $n$, $2\le n\le 4$.

The set of all $L$-domains is partitioned into classes of unimodularly
equivalent domains. For $n \le 3$, there is only one class, i.e. only
one generic $L$-type. For $n=4$, there are 3 generic $L$-types.

The $L$-domain ${\cal A}_n$ of the lattice $A_n^*$ which is dual
to the root lattice $A_n$ is well known. ${\cal A}_n$ has dimension $N$,
i.e. it is generic. It is the unique $L$-domain such that all the
extreme rays of its closure span forms of rank 1. All these facts were
known to Voronoi. He called the domain ${\cal A}_n$ {\em the first type
domain} and one of the forms of $A_n^*$ {\em the principal form of
the first type}.

The $L$-domain of the lattice $A_n$, $n \ge 2$, is a simplicial 
$(n+1)$-dimensional cone; it is described in \cite{BG}. So, nrd$A_n=n+1$ 
for $n \ge 2$. Also nrd$Z_n=n$ for $n \ge 1$ and $A_1, D_2, D_3$ are 
scalings of $Z_1, Z_2, A_3$, respectively.
So, nrd$(A_1=A_1^*)=1$, nrd$(D_2=D_2^*)=2$ and nrd$D_3=4$, nrd$D_3^*=6$.

It is proved in \cite{BG} that the lattice $D_n$ is rigid for $n\ge 4$.

The lattices $E_6$, $E_6^*$, $E_7$, $E_7^*$ and $E_8=E_8^*$ are rigid, 
i.e. their $L$-domains are one-dimensional. The rigidity of root lattices 
$E_6$, $E_7$ and $E_8$ is shown in \cite{BG}. The rigidity of the lattice 
$E_6^*$ was proved independently by Engel and Erdahl (personal 
communications). The rigidity of $E_7^*$ can be proved easily by using a 
nice symmetric quadratic form of $\Gamma({\cal A}^7)=E_7^*$ given in 
\cite{Ba} (see formula (2) there).  

The rigid lattices $A_1$, $E_6$ and $E_7$ are first instances of {\em 
strongly rigid} lattices, i.e. such that amongst of their
Delaunay polytopes there are {\em extreme} ones (see \cite{DGL}), 
i.e. such that any, distinct from a scaling, its affine transformation 
is not a Delaunay polytope. The $1$-simplex and unique Delaunay polytope
of $E_6$ are only such polytopes of dimension at most $6$ (see \cite{DD}).
$A_1$, $D_4=D_4^*$, $E_6$ and $E_7^*$ are rigid lattices, having unique 
type of Delaunay polytope, but only for $A_1$ and $E_6$ this polytope is 
extreme. In \cite{DGL} (see also Chapter 16 of \cite{DL}) were given 10 
examples of extreme Delaunay polytopes: by one of dimension 1, 6, 7, 22. 
23 and by 3 and 2 of dimension 15 and 16, respectively. We believe that 
there is an infinity of them. Moreover, in \cite{DGL} (see also Chapter 
15 of \cite{DL}) was considered the general notion of {\em rank} of a
Delaunay polytope, i.e. the number of degrees of freedom, that it
has, when one deforms it affinely, so that the result remains an 
Delaunay polytope. Clearly, the maximal rank is (as well as maximal nrd) 
$n+1 \choose 2$ and it is realized only by $n$-simplices.

In this note we describe explicitly the $L$-domain for the lattice
$D_n^*$ which is dual to the root lattice $D_n$. This $L$-domain is
special and has dimension $n$ for odd $n$ and dimension 1 for even $n$.
This special $L$-domain is a facet of the closure of several generic 
$L$-domains.

This work completes computation of non-rigidity degree of root lattices and 
their duals. In a sense, it is an addition to the work \cite{CS}, where 
Delaunay and Voronoi polytopes of the root lattices and their duals are 
enumerated. We present the values of nrd for root lattices $L$ and their 
duals in the table below. 

\[\begin{array}{|c||c|c|c|c|c|c|c|c|c|c|c|} \hline
L &A_1=A_1^*& A_n& A_n^*          &D_n   & D_{2m+1}^*&D_{2m}^*
                 & E_6& E_6^* & E_7 & E_7^* & E_8=E_8^* \\ 
               & &n\ge 2& n\ge 1         &n\ge 4& m\ge 2    &m\ge 2
                 &    &       &     &        &         \\  \hline
\mbox{nrd}L    &1& n+1 &\frac{n(n+1)}{2} & 1 & 2m+1 & 1
           &  1 & 1     & 1   &  1     & 1           \\  \hline
 \end{array} \] 
                                               
\section{The cone ${\cal G}_n$}
Let $\{e_i:i \in I_n\}$ be a set of mutually orthogonal vectors of
norms (i.e. of squared lengths) $e_i^2=2\gamma_i$, where $I_n=\{1,2,...,n\}$.
For $S \subseteq I_n$, let $e(S)=\sum_{i \in S}e_i$ and
$\gamma(S)=\sum_{i \in S}\gamma_i$. We introduce the
vector $b$ of norm $\alpha$ as follows:
\begin{equation}
\label{gama}
b=\frac{1}{2}\sum_{i\in I_n}e_i=\frac{1}{2}e(I_n), \mbox{ where }
b^2=\alpha=\frac{1}{2}\sum_{i\in I_n} \gamma_i=\frac{1}{2}\gamma(I_n).
\end{equation}

Let $\overline \gamma$ be the vector with the coordinates
$\{ \gamma_i:1 \le i \le n \}$. Consider the lattice 
$L({\overline \gamma})$ generated by the vector $b$ and any $n-1$ 
vectors $e_i$. If $\gamma_i=1$ for all $i$, and $n\ge 4$, then 
$L({\overline \gamma})=D_n^*$.

We take as a basis of $L({\overline \gamma})$ the vector $b$ and the 
vectors $e_i$ for $1 \le i \le n-1$. Then the coefficients of the 
quadratic form $f_{\overline \gamma}$ corresponding to this basis are 
as follows:
\begin{equation}
     \label{cof}
a_{ii}=e_i^2=2\gamma_i, 1\le i \le n-1,
a_{ij}=e_ie_j=0, 1 \le i,j \le n-1, i \not=j,
\end{equation}
\begin{equation}
\label{con}
a_{nn}=b^2=\alpha, a_{in}=e_i b=\gamma_i, 1\le i \le n-1.
\end{equation}
This form has the following explicit expression
\begin{equation}
\label{fx}
f_{\overline \gamma}(x)=(x_n b+\sum_{i=1}^{n-1} x_ie_i)^2=\alpha x_n^2+
2\sum_1^{n-1}\gamma_ix^2_i+2\sum_1^{n-1}\gamma_ix_ix_n.
\end{equation}
In the basis $\{e_i: i \in I_n\}$, each vertex of
$L({\overline \gamma})$ has integer or half-integer coordinates.

Let $n$ be odd, say $n=2m+1$. Suppose that the parameters $\gamma_i$
satisfy the following $n \choose m$ inequalities
\begin{equation}
\label{dom}
\sum_{i \in S}\gamma_i < \alpha, \mbox{  } S\subset I_n, \mbox{ }
|S|=m.
\end{equation}
Denote by ${\cal G}_n$ the $n$-dimensional domain determined in the
space of variables $\gamma_i$, $i \in I_n$, by the inequalities
(\ref{dom}). Since these inequalities are linear and homogeneous
(recall that $\alpha=\frac{1}{2}\gamma(I_n)$), ${\cal G}_n$ is an
open polyhedral cone. Since $\gamma(S)+\gamma(I_n-S)=2\alpha$,
the inequalities (\ref{dom}) imply the following inequalities
\begin{equation}
\label{mod}
\gamma(T)>\alpha, \mbox{  }T\subset I_n, \mbox{  }|T|=m+1.
\end{equation}
For a set $T$ of cardinality $|T|=m+1$, let $T=S\cup\{i\}$, where
$|S|=m$. Then (\ref{dom}) and (\ref{mod}) imply
\[\alpha<\gamma(T)=\gamma(S)+\gamma_i, \mbox{ i.e. }
\gamma_i>\alpha-\gamma(S)>0. \]
Hence the cone ${\cal G}_n$ lies in the positive orthant of ${\bf R}^n$.

Consider the closure cl${\cal G}_n$ of the cone ${\cal G}_n$. Obviously, 
cl${\cal G}_n$ is defined by the non-strict version of inequalities
(\ref{dom}). So, using that $2\alpha=\gamma(I_n)$, we have
\begin{equation}
\label{clg}
{\rm cl}{\cal G}_n=\{{\overline \gamma}:\gamma(S)-\gamma(I_n-S)\le 0,
\mbox{  }S\subset I_n, |S|=m \}.
\end{equation}
Note that the zero vector belongs to ${\rm cl}{\cal G}_n$.
The automorphism group of cl${\cal G}_n$ is isomorphic to the group
of all permutations of the set $I_n$.

Obviously, the hyperplanes supporting facets of cl${\cal G}_n$ are
contained among the hy\-per\-pla\-nes defined by the equalities
\begin{equation}
\label{fct}
\gamma(S)=\gamma(I_n-S)=\alpha=\frac{1}{2}\gamma(I_n),
\mbox{  }S\subset I_n, \mbox{  }|S|=m.
\end{equation}
Note that the equality $\gamma(S_1)=\alpha$ can be transposed into
the equality $\gamma(S_2)=\alpha$ by the automorphism group, for any
$S_1,S_2 \subset I_n$ with $|S_1|=|S_2|=m$. Hence each of the
equations of (\ref{fct}) determines a facet of cl${\cal G}_n$.

\begin{prop}
\label{extr}
Let $n$ be odd and $n=2m+1\ge 5$, i.e. $m\ge 2$. Then the closure of
${\cal G}_n$ has the following $2n$ extreme rays
\[{\overline \gamma}_q^k=\{\gamma_i=\gamma\ge 0, i \in I_n-\{k\},
\gamma_k=2q\gamma \}, \mbox{  }q=0,1, \mbox{  }k \in I_n. \]
\end{prop}
{\bf Proof}. Let ${\overline \gamma} \in {\rm cl}{\cal G}_n$ be fixed.
Then ${\overline \gamma}=(\gamma_1, \gamma_2,...,\gamma_n)$ defines a
partition of the set $I_n$ as follows. Let the coordinates $\gamma_i$
take $k$ distinct values $0 \le \beta_1<\beta_2<...<\beta_k$, where $k$
is an integer between 1 and $n$. For $1 \le j \le k$, set
$S_j=\{i \in I_n: \gamma_i=\beta_j\}$ and $s_j=|S_j|$. Then
$\sum_{j=1}^ks_j=n=2m+1$ and $I_n=\cup_{j=1}^kS_j$ is the above mentioned 
partition.

Consider the values of $\gamma(S)$ for $S \subset I_n$, $|S|=m$.
$\gamma(S)$ takes a maximal value for the following sets $S$. Let $j_0$
be such that $\sum_{j=1}^{j_0-1}s_j<m+1$, but
$\sum_{j=1}^{j_0}s_j\ge m+1$. Then $\sum_{j=j_0+1}^k s_j \le m$. Let
$S_{max}(T)=T\cup_{j=j_0+1}^k S_j$, where $T \subseteq S_{j_0}$,
$|T|=t_0$ and $t_0:=m-\sum_{j=j_0+1}^k s_j$.
Obviously, $\gamma(S)$ takes the maximal value
$t_0\beta_{j_0}+\sum_{j=j_0+1}^k s_j\beta_j$ if $S=S_{max}(T)$ for
any $T \subseteq S_{j_0}$ of cardinality $|T|=t_0$.

For given ${\overline \gamma} \in {\rm cl}{\cal G}_n$, let
${\cal S}_{max}({\overline \gamma})$ be the system of equations of type
(\ref{fct}), where $S=S_{max}(T)$ for all $T\subseteq S_{j_0}$ with
$|T|=t_0$. If $\overline \gamma$ is an extreme ray of cl${\cal G}_n$, then
${\cal S}_{max}({\overline \gamma})$ determines uniquely up to a multiple
the vector $\overline \gamma$. It is not difficult to see that
${\cal S}_{max}({\overline \gamma})$ can uniquely determine 
$\overline \gamma$ only if $k=2$. Consider this case in detail.

If $k=2$, we have $I_n=S_1\cup S_2$ and $n=2m+1=s_1+s_2$. Let
$s_2\le m$, i.e. $j_0=1$. Then ${\cal S}_{max}({\overline \gamma})$ 
consists of the following equations
\[\gamma(T\cup S_2)=\gamma(T)+\gamma(S_2)=\alpha=\gamma(S_1-T),
\mbox{  }T\subseteq S_1, \mbox{  }|T|=m-s_2. \]
Note that $\gamma_i$ for $i \in S_2$ belongs to the above system only
as a member of the sum $\gamma(S_2)=\sum_{i\in S_2}\gamma_i$. Hence
such a system can determine the coordinates $\gamma_i$, $i \in S_2$,
only if $s_2=1$.

Now the above system implies that $\gamma_i$ takes the same value, say
$\gamma$, for all $i \in S_1$. In fact, let $i_1,i_2 \in S_1$,
$i_1 \in T_1$, $i_1 \not \in T_2$, $i_2 \not \in T_1$, $i_2 \in T_2$,
for some $T_1,T_2 \subset S_1$ of cardinality $m-1$. Such $T_1$ and
$T_2$ exist, since $|T_j|=m-1\ge 1$. Subtracting the equation of the
above system for $T=T_2$ from the equation for $T=T_1$, we obtain the
equality $\gamma_{i_1}-\gamma_{i_2}=\gamma_{i_2}-\gamma_{i_1}$, i.e.
$\gamma_{i_1}=\gamma_{i_2}$.

In this case, the above system, where $S_2=\{k\}$, gives
$\gamma_k=\gamma(S_1)-2\gamma(T)=s_1\gamma-2(m-1)\gamma=2\gamma$.
We obtain the extreme ray ${\overline \gamma}_1^k$.

Now, let $s_2>m$, i.e. $s_1<m+1$ and $j_0=2$. A similar analysis shows
that $s_1=1$, say $S_1=\{k\}$, and $\gamma_i$ take the same value,
say $\gamma$, for all $i\in S_2$. This gives $\gamma_k=0$, and we obtain
the extreme ray ${\overline \gamma}_0^k$. The result follows. 

\vspace{3mm}
The facet defined by the equation $\gamma(S)=\alpha$, $|S|=m$, contains
the following $n=2m+1$ extreme rays: ${\overline \gamma}_0^k$,
$k \not \in S$, ${\overline \gamma}_1^k$, $k \in S$. Each facet has the
following geometrical description. The $m+1$ rays ${\overline\gamma}^k_0$,
$k \not \in S$, form an $(m+1)$-dimensional simplicial cone. Similarly,
the $m$ rays ${\overline\gamma}_1^k$, $k \in S$, form an $m$-dimensional
cone. Both these cones intersect by the ray
$\{{\overline\gamma}:\gamma_i=(m+1)\gamma, i \in S, \gamma_i=m\gamma, 
i \not \in S\}$. Hence the cone
hull of these two cones is a cone of dimension $(m+1)+m-1=2m$. This
cone is just a facet of ${\cal G}_n$ for $n=2m+1$.

Let $n$ be even, $n=2m$. In this case, all the inequalities (\ref{dom})
imply the following set of inequalities
\[\gamma(I_n-S)=\gamma(T)>\alpha, \mbox{  }T=I_n-S\subset I_n, |T|=m. \]
We see that this system of inequalities contradicts to the system
(\ref{dom}). This means that the open cone ${\cal G}_n$ for even $n=2m$
is empty. But the solution of the set of equalities (\ref{fct}) is
not empty. Namely, it has the solution $\gamma_i=\gamma \ge 0$ for
all $i \in I_n$. In other words, cl${\cal G}_n$ is the following ray
\[{\rm cl}{\cal G}_{2m}=\{{\overline \gamma}:\gamma_i=\gamma \ge 0,
i \in I_{2m} \}. \]

\section{The domain ${\cal D}_n$}
Denote by ${\cal D}_n$ the domain of forms $f_{\overline\gamma}$,
where $\overline\gamma$ belongs to ${\cal G}_n$.

We prove the following theorem.
\begin{theor}
\label{main}
Let $n$ be odd, $n=2m+1$. The domain ${\cal D}_n$ is an $L$-domain.
It lies in an $n$-dimensional space which is an intersection of
$n \choose 2$ hyperplanes given by the following equalities
\begin{equation}
\label{ij}
a_{ij}=0, 1\le i<j \le n-1, 2a_{in}=a_{ii}, 1 \le i \le n-1.
\end{equation}
The domain ${\cal D}_n$ is cut from this space by the following
inequalities
\begin{equation}
\label{ii}
\sum_{i \in S}a_{ii}<2a_{nn}, \mbox{  }S \subset I_{n-1}, |S|=m,
\end{equation}
\begin{equation}
\label{kk}
2a_{nn}<\sum_{i \in T}a_{ii}, \mbox{ }T \subset I_{n-1}, |T|=m+1.
\end{equation}
There is a one-to-one correspondence between ${\cal D}_n$ and the
cone ${\cal G}_n$ given by the equalities (\ref{cof}) and (\ref{con}).

In particular, the closure of ${\cal D}_n$ has $2n$ extreme rays
$f_0^k$, $f_1^k$, $k \in I_n$, with the coefficients $a_{ij}$ of
these forms defined as follows (where the term $a_{kk}(f_{0,1}^k)$
should be omitted if $k=n$):
\[a_{ii}(f_0^k)=2\gamma, i \in I_{n-1}, i \not=k, a_{kk}(f_0^k)=0,
a_{nn}(f_0^k)=m\gamma; \]
\[a_{ii}(f_1^k)=2\gamma, i \in I_{n-1}, i \not=k, a_{kk}(f_1^k)=4\gamma,
a_{nn}(f_1^k)=(m+1)\gamma; \]
$a_{ij}(f_{0,1}^k)$ for $i\not=j$ are defined by the equations
(\ref{ij}).

The inequalities (\ref{ii}) and (\ref{kk}) define facets of the
closure cl${\cal D}_n$. All facets are domains of
equivalent  $L$-types, each having $n$ extreme rays $f^k_1$, $k \in S$,
$f^k_0$, $k \not \in S$, $S \subset I_{n-1}$, $|S|=m$, or $S=I_n-T$ and
$T$ is as in (\ref{kk}).

If $n$ is even, $n=2m$, then cl${\cal D}_n$ is one dimensional. The ray
cl${\cal D}_{2m}$ is the intersection of the $n \choose 2$ hyperplanes
(\ref{ij}) and the $n-1$ hyperplanes given by the following equalities
\begin{equation}
\label{ni}
2a_{nn}=ma_{ii}, \mbox{  }1 \le i \le n-1.
\end{equation}
\end{theor}

{\bf Proof} of Theorem~\ref{main} will be proceeded as follows. For a
function $f_{\overline\gamma}$ given by (\ref{fx}), we find the Voronoi
polytope. Take attention that the inequalities (\ref{ii}) and (\ref{kk})
in terms of the parameters $\alpha$ and $\gamma_i$ take the form
(\ref{dom}) for $n \not \in S$ and $n\in S$, respectively. We show that
the face poset of the Voronoi polytope does not change if the parameters
of $f_{\overline\gamma}$ change such that they satisfy (\ref{dom}).

On the other hand, we show if at least one of inequalities (\ref{dom})
holds as equality for parameters of a function $f_{\overline\gamma}$,
then the $L$-type of $f_{\overline\gamma}$ differs from the $L$-type
of $f_{\overline\gamma}\in {\cal D}_n$. This will mean that
${\cal D}_n$ is an $L$-domain.

For to find the Voronoi polytope of $f_{\overline \gamma}$ given by
(\ref{fx}), consider the cosets of $2L$ in the lattice
$L=L({\overline\gamma})$. Let $v=x_nb+\sum_{i \in I_{n-1}}x_ie_i$ be
a vector of $L({\overline\gamma})$. Then this vector belongs to the
coset $Q(S,z)$, where $S \subseteq I_{n-1}$ is the set of indices of
odd coordinates $x_i$ and the number $z \in \{0,1\}$ indicates the
parity of the $b$-coordinate $x_n$ of the vector $v$. Note that the 
vector $e(I_n)=2b$ belongs to the trivial coset $Q(\emptyset,0)=2L$. 
Hence the vectors $e(S)$ and $e(I_n-S)$ belong to the same coset for 
any $S \subseteq I_n$. This coset is $Q(S,0)$ if $n \not \in S$, and
$Q(I_n-S,0)$ if $n \in S$.
In particular, $e_n$ belongs to $Q(I_{n-1},0)$, and it is minimal in
this coset. Moreover, we have $b-e(S)=-(b-e(I_n-S))$. So the
$2^n$ vectors $b-e(S)$, $S \subseteq I_n$, are partitioned into
$2^{n-1}$ pairs of opposite vectors.

Note that $e^2(S)=\sum_{i \in S}e^2_i=2\gamma(S)$ and,
according to (\ref{gama}), $\gamma(S)+\gamma(I_n-S)=2\alpha$. Recall
that ${\cal D}_n$ is the domain of $f_{\overline \gamma}$, where
$\overline\gamma$ belongs to ${\cal G}_n$. Hence $\gamma_i$, $i\in I_n$,
satisfy (\ref{dom}). Taking in attention (\ref{dom}), we see that,
for $|S|\le m$, the norm of $e(S)$ is less than the norm of
$e(I_n-S)=e(T)$ for $|T|\ge m+1$, $S,T \subset I_n$.

If $f_{\overline \gamma}$ go to the boundary of ${\cal D}_n$, then
the sets of minimal vectors of some cosets change. At first we describe
the simple cosets which are constant on the closure of ${\cal D}_n$.
The norm of minimal vectors of a coset is called also norm of the
coset. These are the following cosets:

\vspace{0.3cm}
The $n$ cosets $Q(\{i\},0)$, $i \in I_{n-1}$, and $Q(I_{n-1},0)$ of
norms $2\gamma_i$, $i \in I_n$, with minimal vectors
$e_i$, $i \in I_{n-1}$ and $e_n=2b-e(I_{n-1})$, respectively.

The $2^{n-1}$ cosets $Q(S,1)$ of norm $\alpha$ with minimal vectors
$b-e(S)$, $S \subseteq I_{n-1}$.

\vspace{0.3cm}
The $2^{n-1}$-$n$   non-simple   cosets  $Q(S,0)$,  $S  \subset  I_n$,
$1<|S|\le m$, have norms $\gamma(S)$ with minimal vectors
$\sum_{i \in S}\varepsilon_ie_i$, where $\varepsilon_i\in\{\pm 1\}$.
If $\alpha=\gamma(S)$, then $|S|=m$ and the coset $Q(S,0)$ contains
also the vector $\sum_{i \in I_n-S}\varepsilon_ie_i$.

\vspace{0.3cm}
Recall that the minimal vectors of simple cosets determine facets of the
Voronoi polytope.
Consider a point $x \in {\bf R}^n$ in the basis $\{e_i:i \in I_n\}$,
$x=\sum_{i \in I_n}x_ie_i$. Then $x$ belongs to the Voronoi polytope
$P$ of $L({\overline\gamma})$ if the inequalities
\[-\frac{v^2}{2}\le xv \le \frac{v^2}{2} \]
hold for all minimal vectors $v$ of simple cosets of
$L({\overline\gamma})$. Using (\ref{cof}), (\ref{con}) and the identity
$b=\frac{1}{2}\sum_{i \in I_n}e_i$, we obtain the following system of
inequalities describing the Voronoi polytope of $L({\overline\gamma})$:
\begin{equation}
\label{xv}
-\frac{1}{2}\le x_i \le \frac{1}{2}, i \in I_n,
\end{equation}
\begin{equation}
\label{xg}
-\frac{1}{2}\alpha \le \sum_{i \in I_n}\gamma_i\varepsilon_ix_i
\le \frac{1}{2}\alpha, \mbox{  }
\varepsilon_i \in \{\pm 1\}, i \in I_n.
\end{equation}
Here the inequality (\ref{xg}) is given by minimal vectors of $Q(S,1)$
such that $\varepsilon_i=-1$ if $i \in S$, and $\varepsilon_i=1$ if
$i \not \in S$.

Note that $\sum_{i \in I_n}\gamma_i\varepsilon_i x_i$ is the linear
function on $\varepsilon_ix_i$ taking maximal value if
$\varepsilon_ix_i \ge 0$ for all $i \in I_n$. Hence the right hand
inequality in (\ref{xg}) holds as equality for a vertex $x$ only if
$\varepsilon_ix_i>0$ for $x_i \not=0$.

An analysis of the system (\ref{xv}), (\ref{xg}) shows, that for each
vertex $x$ there is the opposite vertex $-x$, and $x$ has the following
coordinates
\begin{equation}
\label{vert}
x_i=\frac{1}{2}\varepsilon_i, i \in S \subseteq I_n, |S|=m,
x_k=\frac{\varepsilon_k}{2\gamma_k}(\alpha-\gamma(S)),
x_l=0, \mbox{ for } l\in I_n-(S\cup\{k\}).
\end{equation}
There are $n \choose m$ positive vertices of this type. Taking in
attention signs, we obtain $2^{m+1}{n\choose m}$ vertices of the Voronoi
polytope. Denote the vertex (\ref{vert}) by $x(k;S)$.

The form of the vertex $x(k;S)$ shows that some vertices can be glued
if and only if the equality $\alpha=\gamma(S)$ holds for some set $S$.
If $\alpha=\gamma(S)$, then $x_k(k;S)=0$, and the $m+1$ vertices
$x(l;S)$, $l \in I_n-S$, are glued into one vertex.
This means that if $\alpha=\gamma(S)$ for some $S$, then $L$-type of
$f_{\overline\gamma}$ changes. So, we proved that the inequalities
(\ref{dom}), i.e. the inequalities (\ref{ii}) and (\ref{kk}) hold for
$f \in {\cal D}_n$.

Now, we show that ${\cal D}_n$ lies in the intersection of the
hyperplanes (\ref{ij}). It is proved in \cite{BG} that the equations
of the hyperplanes in the intersection of which an $L$-type domain
lies are given by some linear forms on norms of minimal vectors of
cosets of $2L$ in $L$. Some of such linear forms are obtained by
equating norms of minimal vectors of a non-simple coset. There
are $L$-type domains for which linear forms of last type are sufficient
for to describe the space, where this $L$-type domain lies. This is
so in our case.

In fact, it is sufficient to consider the non-simple cosets $Q(S,0)$
for $|S|=2$, i.e. to equate the norms of vectors $e_i+e_j$ and
$e_i-e_j$, $i,j \in I_n$. The equality $(e_i+e_j)^2=(e_i-e_j)^2$ 
implies $e_ie_j=0$, i.e. $a_{ij}=0$, $0 \le i<j \le n-1$. We obtain the 
first equalities in (\ref{ij}).
For $j=n$, we have $e_ie_n=-e_i(2b-e(I_{n-1}))=0$. Since $e_ie_j=0$,
this equality is equivalent to $2be_i=e_i^2$. We obtain the second
equalities in (\ref{ij}). If we set $be_i=\gamma_i$, $b^2=\alpha$,
we obtain the original function $f_{\overline\gamma}$.

Now let $n=2m$ be even. In this case, for $f \in {\rm cl}{\cal D}_n$,
all cosets of $2L$ in $L$ (excluding the cosets $Q(S,0)$ for $|S|=m$)
are the same as in the odd case. But, for $|S|=m$, $Q(S,0)$ contains
beside the vectors $\sum_{i\in S}\varepsilon_i e_i$ also the vectors
$\sum_{i\in I_n-S}\varepsilon_i e_i$. Since norms of these vectors are
$\gamma(S)$ and $\gamma(I_n-S)$, respectively, the equating of these
norms gives the system (\ref{fct}). This system has the unique solution
$\gamma_i=\gamma$ for all $i \in I_n$.

Hence, for $n=2m$, $\alpha=\frac{1}{2}\gamma(I_n)=\frac{1}{2}(\gamma(S)
+\gamma(I_n-S))=\gamma(S)=m\gamma$. Taking in attention (\ref{cof})
and (\ref{con}), we can rewrite this equality as $2a_{nn}=ma_{ii}$ for
any $i\in I_n$. So, we obtain (\ref{ni}). This means that cl${\cal D}_n$
is a ray, which lies in the intersection of the hyperplanes given by
the equations (\ref{ii}), (\ref{kk}) and (\ref{ni}). Any
$f \in {\rm cl}{\cal D}_n$ is a rigid (i.e. edge) form.

So, Theorem~\ref{main} is proved. 

\vspace{3mm}
Recall that $L({\overline\gamma})=D_n^*$ if $\gamma_i=1$ for all
$i \in I_n$. Since, for this $\overline\gamma$, the parameters of
$f_{\gamma}$ satisfy (\ref{dom}), where $\alpha=m+\frac{1}{2}$, this 
implies the following

\vspace{3mm}
{\bf Corollary} {\em ${\cal D}_n$ is the $L$-domain of $D_n^*$.
In particular, the lattice $D_{2m}^*$, $m \ge 2$, is rigid.}

\vspace{3mm}
{\bf Remark}. Note that, for $n$ odd, $n=2m+1$, the extreme rays $f_0^k$ 
have rank $n-1=2m$. These forms are forms of lattices isomorphic to
$\gamma D^*_{n-1}=\gamma D^*_{2m}$.

\vspace{3mm}
We saw in the proof of Theorem 1 that the Voronoi polytope of the
lattice $D_n^*$ is an $n$-cube whose vertices are cut by hyperplanes
(\ref{xg}). A description of the Voronoi polytope of $D_n^*$ can
be found in \cite{CS}.

Theorem 1 is a generalization of the result of \cite{EG}, where the
$L$-domain of the lattice $D_5^*$ is described in detail.


\end{document}